# Induced Measures on $\mu^{**}$-measurable Sets

*Abstract* We investigate extension of a measure to a very general set of undetermined structure. Structure may be imposed on this set in special cases.

## 1. INTRODUCTION

Interest in the extension has continued into recent times with attention being paid to special types of measures, or to extension into particular types of sets, rather than just the general theory of extension of measures. [2] showed that the Lebesgue-like extension of every finitary measure on the Cartesian product of a countable number of discrete topological spaces are a measure on the lattice of open sets. [6] Considered *G*-invariant measures, where *G* is and at most countable group of bijections. [3] Established that each sequentially continuous normed measure on a bold algebra of fuzzy sets can be uniquely extended to sequentially continuous measure on the generated Lukasiewicz tribe. [4] Also characterized extension of probability measures as a completely categorical construction. We consider extension of a measure onto a very general set, so there is great variety in its possible structure. In a particular situation, appropriate structure can be imposed.

A measure $\mu$ on a ring *R* induces an outer measure $\mu^*$ on the smallest $\sigma$-set $H(R)$ containing *R*. $\mu^*$ in turn induces a complete measure on $\overline{S}$, the $\sigma$-algebra of all $\mu^*$-measurable sets (c.f.[5]). We investigate whether $\mu$ may induce a measure on a very general set. We show that $\mu$ induces a complete measure on a set containing $\overline{S}$ which arises in a natural way by considering sets whose intersections with $\mu^*$-measurable sets are also $\mu^*$-measurable. This was suggested by the definition of sets measurable with respect to a measurable space (c.f. [1]). The results obtained are generalizations of results from [5].

## 2. DEFINITIONS AND NOTATIONS

<u>Definition</u>  Let $(G, \mathbf{B})$ be a measurable space. Let $A \subset G$, and suppose that for every $B \in \mathbf{B}$, we have $A \cap B \in \mathbf{B}$. We then say that $A$ is <u>measurable</u> with respect to $(G, \mathbf{B})$. We denote the totality of such measurable sets $A$ by $\widetilde{\mathbf{B}}$. (c.f. [1]).

The following definition and two notations are from [5].

<u>Notation</u>  For any class $\mathbf{E}$ of sets, $\mathbf{S}(\mathbf{E})$ is the smallest $\sigma$-ring containing $\mathbf{E}$.

<u>Definition</u>  A non-empty class $\mathbf{E}$ of sets is hereditary if

$F \subset \mathbf{E}$, whenever $E \subset \mathbf{E}$ and $F \subset E$.

<u>Notation</u>  $H(R)$ denotes the smallest hereditary set containing a ring $R$.

<u>Notation</u>  Let $X$ be a set having non-empty intersection with the set containing all elements occurring in the sets of $R$.

<u>Notation</u>  Let $K$  $K = \{B \in X | B \cap E \in H(R) \ \forall \ E \in H(R)\}$

<u>Definition</u>  A set $Q \in K$ is $\mu^{**}$-measurable if $Q \cap E$ is $\mu^*$-measurable $\forall \ E \in H \ni E$ is $\mu^*$-Measurable; that is, if $\forall A \in H$ and for any $\mu^*$-measurable set $E$,

$\mu^*(A) = \mu^*[A \cap (Q \cap E)] + \mu^*[A \cap (Q \cap E)^c]$.

## 3. INDUCED MEASURES

**Theorem 3.1** If $\mu^*$ is an outer measure on a hereditary $\sigma$-ring $H$ and if $\overline{\overline{S}}$ is the class of all $\mu^{**}$-measurable sets, then $\overline{\overline{S}}$ is a ring.

**Proof** Let $L, M \in \overline{\overline{S}}$. Then for any $G \in \overline{S}$, $L \cap G \in \overline{S}$ and $M \cap G \in \overline{S}$.

$$\therefore (L \cup M) \cap G \in \overline{S}$$

$$\therefore L \cup M \in \overline{\overline{S}} \qquad (I)$$

We now show that $(M - L) \cap G \in \overline{S}, \forall\, G \in \overline{S}$.

Since $M \cap G$ is $\mu^*-$ measurable,

$$\mu^*(A) = \mu^*[A \cap (M \cap G)] + \mu^*[A \cap (M \cap G)'] \qquad (II)$$

Since $L \cap G$ is $\mu^*-$ measurable,

$$\mu^*[A \cap (M \cap G)] = \mu^*[A \cap (M \cap G) \cap (L \cap G)] + \mu^*[A \cap (M \cap G) \cap (L \cap G)'] \quad (III)$$

Also, $\mu^*[A \cap (M \cap G)']$

$$= \mu^*[(A \cap (M \cap G)' \cap (L \cap G)] + \mu^*[A \cap (M \cap G)' \cap (L \cap G)'] \qquad [IV]$$

Substituting (III) and (IV) into (II)

$$\mu^*(A)$$

$$= \mu^*[A \cap (M \cap G) \cap (L \cap G)]$$

$$+ \mu^*[A\cap(M\cap G)\cap(L\cap G)']$$

$$+ \mu^*[A\cap(M\cap G)'\cap(L\cap G)]$$

$$+ \mu^*[A\cap(M\cap G)'\cap(L\cap G)'] \tag{V}$$

It is not difficult to establish that:

$$A\cap(M\cap G)\cap(L\cap G)' = A\cap[(M-L)\cap G] \tag{VI}$$

Also, substituting $A\cap[(M\cap L')\cap G]'$ in place of $A$ in (V).

$$\mu^*\{A\cap[(M-L)\cap G]'\}$$

$$= \mu^*\{A\cap[(M\cap L')\cap G]'\}$$

$$= \mu^*\{A\cap[(M\cap L')\cap G]' \cap (M\cap G)\cap(L\cap G)\}$$

$$+ \mu^*\{A\cap[(M\cap L')\cap G]' \cap [(M\cap G)\cap(L\cap G)']\}$$

$$+ \mu^*\{A\cap[(M\cap L')\cap G]' \cap (M\cap G)'\cap(L\cap G)\}$$

$$+ \mu^*\{A\cap[(M\cap L')\cap G]' \cap (M\cap G)'\cap(L\cap G)'\} \tag{VII}$$

We can easily prove (1), (2), (3) and (4) below.

$$A\cap[(M\cap L')\cap G]' \cap (M\cap G)\cap(L\cap G) = A\cap(M\cap G)\cap(L\cap G) \tag{1}$$

$$A \cap [(M \cap L') \cap G]' \cap (M \cap G) \cap (L \cap G)' = \phi \qquad (2)$$

$$A \cap [(M \cap L') \cap G]' \cap (M \cap G)' \cap (L \cap G) = A \cap (M \cap G)' \cap (L \cap G) \qquad (3)$$

$$A \cap [(M \cup L') \cap G]' \cap (M \cap G)' \cap (L \cap G)' = A \cap (M \cap G)' \cap (L \cap G)' \qquad (4)$$

From (V), (VI), (1), (2), (3), and (4)

$$\mu^*(A) = \mu^*\{A \cap [(M-L) \cap G]\} + \mu^*\{A \cap [(M-L) \cap G]'\}$$

Hence $(M - L) \cap G$ is $\mu^*-$ measurable for any $G \in \overline{S}$.

Hence $(M-L)$ is $\mu^{**}-$ measurable i.e. $(M\text{-}L) \in \overline{\overline{S}}$

(VIII)

From (I) and (VIII) $\overline{\overline{S}}$ is a ring. □

<u>Theorem 3.2</u>  If $\mu^*$ is an outer measure on a hereditary $\sigma -$ ring $H$ and if $\overline{\overline{S}}$ is in the class of all $\mu^{**}-$ measurable sets, then $\overline{\overline{S}}$ is a $\sigma -$ ring.

<u>Proof:</u> Let $L_1, L_2, \ldots$ be an infinite sequence of sets in $\overline{\overline{S}}$ Then for any

$$G \in \overline{S}, \ L_i \cap G \in \overline{S}; i = 1, 2, \ldots$$

$$(L_1 \cup L_2 \cup \ldots) \cap G = (L_1 \cap G) \cup (L_2 \cap G) \cup \ldots \in \overline{S}$$

since $\overline{S}$ is an $r$-ring (Theorem B

Hence $(L_1 \cup L_2 \cup ...) \in \overline{\overline{S}}$

Hence $\overline{\overline{S}}$ is a $\sigma$-ring. □

Definition Let $\tilde{\mu}$ be a set function defined on $\overline{\overline{S}}$ by

$$\tilde{\mu}(P) = \sup_{T \in \overline{S}} \mu^*(P \cap T), \forall P \in \overline{\overline{S}}$$

Lemma 3.3 If $A \in H$ and if $\{L_n\}$ is a disjoint sequence of sets in $\overline{\overline{S}}$ with

$$\bigcup_{n=1}^{\infty} L_n = L, \text{ then}$$

$$\sup_{T \in \overline{S}} \sum_{n=1}^{\infty} \mu^*(A \cap L_n \cap T) = \sum_{n=1}^{\infty} \sup_{T \in \overline{S}} \mu^*(A \cap L_n \cap T)$$

Proof Let $\mu^*(A \cap T_n \cap T)$ take its maximum value for $T = T_n$ and

let $\bigcup_{n=1}^{\infty} T_n = U.$ Then $U \in \overline{S}$ and

$$\sup_{T \in \overline{S}} \mu^*(A \cap L_n \cap T) = \mu^*(A \cap L_n \cap U) \tag{I}$$

Suppose that $\exists V \in \overline{S} \ni \sum_{n=1}^{\infty} \mu^*(A \cap L_n \cap V) > \sum_{n=1}^{\infty} \mu^*(A \cap L_n \cap U).$

Then $\exists$ some value of $n$, $N$ says, $\ni$

$\mu^*(A \cap L_n \cap V) > \mu^*(A \cap L_n \cap U)$, contradicting (I).

Hence $\sup\limits_{T \in \overline{S}} \sum\limits_{n=1}^{\infty} \mu^*(A \cap L_n \cap T) = \sum\limits_{n=1}^{\infty} \mu^*(A \cap L_n \cap U)$

$$= \sum\limits_{n=1}^{\infty} \sup\limits_{T \in \overline{S}} \mu^*(A \cap L_n \cap T) \qquad \square$$

<u>Theorem 3.4</u>  If $A \in H$ and if $\{L_n\}$ is a disjoint sequence of sets in $\overline{\overline{S}}$ with $\bigcup\limits_{n=1}^{\infty} L_n = L$, then

$$\tilde{\mu}(A \cap L) = \sum\limits_{n=1}^{\infty} \tilde{\mu}(A \cap L_n)$$

<u>Proof</u>  Let $T$ be an arbitrary element of $\overline{S}$. Then $\{L_n \cap T\}$ is a disjoint sequence of sets in $\overline{S}$ and $L \cap T = \bigvee\limits_{n=1}^{\infty}(L_n \cap T)$.

Hence

$$\mu^*(A \cap L \cap T) = \sum\limits_{n=1}^{\infty} \mu^*(A \cap L_n \cap T)$$

from Theorem 11.B, by [5]

$$\therefore \sup\limits_{T \in \overline{S}} \sum\limits_{n=1}^{\infty} \mu^*[A \cap (L_n) \cap T]$$

$$= \sum\limits_{n=1}^{\infty} \sup\limits_{T \in \overline{S}} \mu^*(A \cap L_n \cap T) \text{ from the previous lemma.}$$

i.e. $\tilde{\mu}(A \cap L) = \sum_{n=1}^{\infty} \tilde{\mu}(A \cap L_n)$      □

Definition  The set function $\bar{\mu}$ is defined on $\bar{S}$ by

$$\bar{\mu}(E) = \mu^*(E), \text{ for } E \in \bar{S}.$$

Remark  $\bar{\mu}$ is a complete measure on $\bar{S}$. (Theorem 11 (c), [5]).

Lemma 3.5  $\sup\limits_{T \in \bar{S}} \left\{ \sum_{n=1}^{\infty} \bar{\mu}(L_n \cap T) \right\} = \sum_{n=1}^{\infty} \sup\limits_{T \in \bar{S}} \bar{\mu}(L_n \cap T)$

where $\{L_n\}$ is a disjoint sequence in $\bar{\bar{S}}$.

Proof  As $\bar{S}$ is monotone, it is easy to show that $\exists T_n \ni \bar{\mu}(L_n \cap T)$ takes its maximum value for $T = T_n$. Let $U = \bigcup_{n=1}^{\infty} T_n$.

Then $U \in \bar{S}$ and

$$\sup\limits_{T \in \bar{S}} \bar{\mu}(L_n \cap T) = \bar{\mu}(L_n \cap U)$$

$$\therefore \sum_{n=1}^{\infty} \sup\limits_{T \in \bar{S}} \bar{\mu}(L_n \cap T) = \sum_{n=1}^{\infty} \bar{\mu}(L_n \cap U)$$

$$= \bar{\mu}\left[ \left( \bigcup_{n=1}^{\infty} L_n \right) \cap U \right]$$

$$= \sup\limits_{T \in \bar{S}} \bar{\mu}\left[ \left( \bigcup_{n=1}^{\infty} L_n \right) \cap T \right]$$

$$= \sup_{T \in \overline{\overline{S}}} \left\{ \sum_{n=1}^{\infty} \overline{\mu}(L_n \cap T) \right\}$$

□

Theorem 3.6  If $\mu^*$ is an outer measure on a hereditary $\sigma$-ring H and if $\overline{\overline{S}}$ is the class of all $\mu^{**}$-measurable sets, then every set of outer measure zero belongs to $\overline{\overline{S}}$ and $\tilde{\mu}$ is a complete measure on $\overline{\overline{S}}$.

Proof: If $E \in H$ and $\mu^*(E) = 0$, then for any $G \in \overline{S}$ and $A \in H$,

$$\mu^*(A) = \mu^*(E) + \mu^*(A) \geq \mu^*[A \cap (E \cap G)] + \mu^*\left[A \cap (E \cap G)^c\right]$$

Since

$$\mu^*(A) \leq \mu^*[A \cap (E \cap G)] + \mu^*\left[A \cap (E \cap G)^c\right],$$

$$\mu^*(A) = \mu^*[A \cap (E \cap G)] + \mu^*\left[A \cap (E \cap G)^c\right].$$

Hence $E \in \overline{\overline{S}}.$

Countable Additivity.  Let $\{L_n\}$ be a disjoint sequence in $\overline{\overline{S}}$. For any $T \in \overline{S}$,

$$\overline{\mu}\left[\left(\bigcup_{n=1}^{\infty} L_n\right) \cap T\right] = \overline{\mu}\left[\bigcup_{n=1}^{\infty} (L_n \cap T)\right]$$

$$= \sum_{n=1}^{\infty} \overline{\mu}(L_n \cap T)$$

since $\bar{\mu}$ is a complete measure on $\bar{S}$. (Theorem 11.c [5]).

$$\therefore \sup_{T \in \bar{S}} \bar{\mu}\left[\left(\bigcup_{n=1}^{\infty} L_n\right) \cap T\right]$$

$$= \sup_{T \in \bar{S}} \left\{\sum_{n=1}^{\infty} \mu(L_n \cap T)\right\}$$

$$= \sum_{n=1}^{\infty} \sup_{T \in \bar{S}} \bar{\mu}(L_n \cap T)$$

(by the previous lemma).

$$\therefore \sup_{T \in \bar{S}} \mu^*\left[\left(\bigcup_{n=1}^{\infty} L_n\right) \cap T\right]$$

$$= \sum_{n=1}^{\infty} \sup_{T \in \bar{S}} \mu^*(L_n \cap T)$$

i.e. $\tilde{\mu}\left(\bigcup_{n=1}^{\infty} L_N\right) = \sum_{n=1}^{\infty} \tilde{\mu}(L_n)$.

$\therefore \tilde{\mu}$ is count ably additive and hence a measure .

Completeness

If $E \in \bar{\bar{S}}, F \subset C$ and $\tilde{\mu}(E) = 0$, then $\tilde{\mu}(F) = 0$ and so $F \in \bar{\bar{S}}$. Hence $\tilde{\mu}$ is complete. □

Remark $\tilde{\mu}$ is called the measure induced by $\mu^*$.

Theorem 3.7 Every set in $S(R)$ is $\mu^{**}$-measurable.

Proof: Let $E \in \overline{S}$. For any element $G$ of $\overline{S}$, $E \cap G \in \overline{S}$. Hence E is

$\mu^{**}$- measurable.

$$\therefore E \in \overline{\overline{S}}$$

$$\therefore \overline{S} \subset \overline{\overline{S}}$$

Since $R \subset \overline{S}$ (Theorem A, page 49, [5]), $R \subset \overline{\overline{S}}$.

Since $\overline{\overline{S}}$ is a $\sigma$-ring.

$$S(R) \subset \overline{\overline{S}}. \qquad \square$$

Theorem 3.8 If $E \in H(R)$, then $\mu^*(E) = \inf\{\tilde{\mu}(F) : E \subset F \in \overline{S}\}$

$$= \inf\{\tilde{\mu}(F) : E \subset F \in S(R)\}$$

Proof: Recall that

$$\mu^*(F) = \inf\left\{\sum_{n=1}^{\infty} \mu(E_n) : E_n \in R \text{ and } E \subset \bigcup_{n=1}^{\infty} E_n\right\}$$

If $F \in R$, then, by the above definition, $\mu^*(F) = \mu(F)$. Since $F \in R$,

$F \in \overline{S}, \ F \in \overline{\overline{S}}$.

$\therefore \tilde{\mu}(F) = \sup_{T \in \overline{S}} \mu^*(F \cap T) \ = \ \mu^*(F)$

$\therefore$ if $F \in R, \ \mu(F) = \tilde{\mu}(F)$

From (I),

$$\mu^*(E) \geq \inf \left\{ \sum_{n=1}^{\infty} \tilde{\mu}(E_n) : E_n \in S(R), \text{ and } E \subset \bigcup_{n=1}^{\infty} E_n \right\}.$$

Since every sequence $\{E_n\}$ of sets in $S(R)$ for which

$$E \subset \bigcup_{n=1}^{\infty} E_n = F$$

may be replaced by a disjoint sequence with the same property without increasing the sum of the measures of the terms of the sequence, and since, by the definition of $\tilde{\mu}$,

$\tilde{\mu}(E) \geq \mu^*(F), \forall F \in \overline{S}$

it follows that

$\mu^*(E) \geq \inf\{\tilde{\mu}(F) : E \subset F \in S(R)\}$

$\geq \inf\{\tilde{\mu}(F) : E \subset F \in \overline{S}\}$

$\geq \inf\{\mu^*(F) : E \subset F \in \overline{S}\}$

$$\geq \mu^*(E)$$

and the result follows. □

Remark  Given $n$ measure spaces $(X_i, S_i, \mu_i)$, with the $S_i,s$ being mutually disjoint, the measure a $\mu_T$ can be defined on the space T, where

$$T = \left\{\bigcup_{i=1}^n A_i \Big| A_i \in S_i, i = 1,2,...,n.\right\}$$

and $\mu_T\left(\bigcup_{i=1}^n A_i\right) = \sum_{i=1}^n A_i$

$T$ can be defined in a variety of other ways.

We can also consider the outer measure $\mu_i^*$ and set function $\tilde{\mu}_i$ associated with each $\mu_i . i = 1,2,...,n,$, and investigate the measure induced by the $\mu_i's$.

**APPENDIX**

**Detailed Proof of Theorem 3.1**

Proof  Let $L, M \in \overline{\overline{S}}.$ Then for any $G \in \overline{S}, L \cap G$ and $M \cap G \in \overline{S}.$

Now, $(L \cup M) \cap G = (L \cap G) \cup (M \cap G)$

Since $L \cap G, M \cap G \in \overline{S}$ and $\overline{S}$ is a ring.

$(L \cap G) \cup (M \cap G) \in \overline{S}$

$\therefore (L \cup M) \cap G \in \overline{S}$

$\therefore L \cup M \in \overline{\overline{S}}$ \hfill (I)

We now show that $(M - L) \cap G \in \overline{S}, \forall\, G \in \overline{S}.$

Since $M \cap G$ is $\mu^* -$ measurable,

$\mu^*(A) = \mu^*\left[A \cap (M \cap G)\right] + \mu^*\left[A \cap (M \cap G)'\right]$ \hfill (II)

Since $L \cap G$ is $\mu*-$ measurable,

$\mu*[A \cap (M \cap G)]$

$= \mu*[A \cap (M \cap G) \cap (L \cap G)] + \mu*[A \cap (M \cap G) \cap (L \cap G)']$ \hfill (III)

Also, $\mu*[A \cap (M \cap G)']$

$= \mu*[(A \cap (M \cap G)' \cap (L \cap G)] + \mu*[A \cap (M \cap G)' \cap (L \cap G)']$ \hfill [IV]

Substituting (III) and (IV) into (II)

$\mu*(A)$

$= \mu*[A \cap (M \cap G) \cap (L \cap G)]$

$+ \mu*[A \cap (M \cap G) \cap (L \cap G)']$

$+ \mu*[A \cap (M \cap G)' \cap (L \cap G)]$

$+ \mu*[A \cap (M \cap G)' \cap (L \cap G)']$ \hfill (V)

Now $A \cap (M \cap G) \cap (L \cap G)'$

$= A \cap (M \cap G) \cap (L' \cup G')$

$= (A \cap M \cap G \cap L') \cup (A \cap M \cap G \cap G'$

$= A \cap [(M \cap L') \cap G]$

$= A \cap [(M - L) \cap G]$ \hfill (VI)

Also, $\mu^*\{A \cap [(M - L) \cap G]'\}$

$= \mu^*\{A \cap [(M \cap L') \cap G]'\}$

$= \mu^*\{A \cap [(M \cap L') \cap G]' \cap (M \cap G) \cap (L \cap G)\}$

$+ \mu^*\{A \cap [(M \cap L') \cap G]' \cap [(M \cap G) \cap (L \cap G)']\}$

$+ \mu^*\{A \cap [(M \cap L') \cap G]' \cap (M \cap G)' \cap (L \cap G)\}$

$+ \mu^*\{A \cap [(M \cap L') \cap G]' \cap (M \cap G)' \cap (L \cap G)'\}$, (VII)

substituting $A \cap [(M \cap L') \cap G]'$ in place of $A$ in (V).

$A \cap [(M \cap L') \cap G]' \cap (M \cap G) \cap (L \cap G)$

$= A \cap (M' \cup L \cup G') \cap (M \cap G) \cap (L \cap G)$

$= A \cap [(M \cap G)' \cup L] \cap (M \cap G) \cap (L \cap G)$

$= A \cap (M \cap G) \cap (L \cap G)$ (1)

$A \cap [(M \cap L') \cap G]' \cap (M \cap G) \cap (L \cap G)'$

$= A \cap (M' \cup L \cup G') \cap (M \cap G) \cap (L' \cup G')$

$$= A \cap (M' \cup L \cup G') \cap (M \cap G \cap L')$$

$$= A \cap (M \cap G \cap L')' \cap (M \cap G \cap L')$$

$$= \phi \qquad (2)$$

$$A \cap [(M \cap L') \cap G]' \cap (M \cap G)' \cap (L \cap G)$$

$$= A \cap (M' \cup L \cup G') \cap (M' \cup G') \cap (L \cap G)$$

$$= A \cap (M' \cup L \cup G') \cap (M' \cap L \cap G)$$

$$= A \cap (M' \cap L \cap G) \cup (M' \cap L \cap G) \cup \phi$$

$$= A \cap (M' \cap L \cap G')$$

$$= A \cap (M \cap G)' \cap (L \cap G) \qquad (3)$$

$$A \cap [(M \cup L') \cap G]' \cap (M \cap G)' \cap (L \cap G)'$$

$$= A \cap (M' \cup L \cup G') \cap (M' \cup G') \cap (L' \cup G')$$

$$= A \cap (M' \cup L \cup G') \cap \{[M' \cap (L' \cup G')] \cup [G' \cap (L' \cup G')]\}$$

$$= A \cap [(M \cap G)' \cup L] \cap \{[M' \cap (L \cap G)'] \cup [G \cap (L \cap G)']\}$$

$$= A \cap [(M \cap G)' \cup L] \cap [[\{M' \cup [G' \cap (L \cap G)']\} \cap \{(L \cap G)' \cup [G' \cap (L \cap G)']\}]]$$

$$= A \cap [(M \cap G)' \cup L] \cap \left[ (M \cap G)' \cap (L \cap G)' \right]$$

$$= A \cap \left\{ \left[ (M \cap G)' \cap (L \cap G)' \right] \cup \left[ L \cap (M \cap G)' \cap (L \cap G)' \right] \right\}$$

$$= A \cap \left\{ \left[ (M \cap G)' \cap (L \cap G)' \right] \cup \left[ L \cap (M \cup G)' \cap (L' \cup G') \right] \right\}$$

$$= A \cap \left\{ \left[ (M \cap G)' \cap (L \cap G)' \right] \cup \left[ (M \cup G)' \cap (L \cap G') \right] \right\}$$

$$= A \cap \left\{ (M \cap G)' \cap \left[ (L \cap G)' \cup (L \cap G') \right] \right\}$$

$$= A \cap \left\{ (M \cap G)' \cap \left[ (L' \cup G') \cup (L \cap G') \right] \right\}$$

$$= A \cap \left\{ (M \cap G)' \cap \left[ U \cap (L' \cap G') \right] \right\}$$

$$= A \cap (M \cap G)' \cap (L \cap G)' \qquad (4)$$

From (V), (VI), (1), (2), (3), and (4)

$$\mu^*(A) = \mu^* \{ A \cap [(M - L) \cap G] \} + \mu^* \{ A \cap [(M - L) \cap G]' \}$$

Hence $(M - L) \cap G$ is $\mu^*-$ measurable for any $G \in \overline{S}$.

Hence (M-L) is $\mu^{**}-$ measurable i.e. (M-L) $\in \overline{\overline{S}}$

(VIII)

From (I) and (VIII) $\overline{\overline{S}}$ is a ring. ☐